\documentclass[12pt]{article}

\usepackage[margin=1in]{geometry}
\usepackage{amsmath,amsthm,amssymb}
\usepackage{color,graphicx}
\usepackage{times}
\usepackage[colorlinks,linkcolor=black,bookmarksopen,bookmarksnumbered,citecolor=black,urlcolor=black]{hyperref}

\newcommand{\bb}{B\&B } 
\newcommand{\bbns}{B\&B} 
\newcommand{\bnb}{branch-and-bound }
\newcommand{\zad}[1]{\mathbb{Z}^{#1}}
\newcommand{\width}{{\rm width}}
\newcommand{\iwidth}{{\rm iwidth}}
\newcommand{\w}{{\rm w}} 

\newcommand{\R}{ {\mathbb R} }
\newcommand{\Z}{ {\mathbb Z} }
\newcommand{\beq}{\begin{equation}}
\newcommand{\eeq}{\end{equation}}
\newcommand{\beqast}{\begin{eqnarray*}}
\newcommand{\eeqast}{\end{eqnarray*}}
\newcommand{\bdef}{\begin{definition}}
\newcommand{\Edef}{\end{definition}}
\newcommand{\ba}{\begin{array}}
\newcommand{\ea}{\end{array}}
\newcommand{\bex}{\begin{example}}
\newcommand{\eex}{\end{example}}
\newcommand{\ble}{\begin{lemma}}
\newcommand{\ele}{\end{lemma}}
\newcommand{\bth}{\begin{theorem}}
\newcommand{\enth}{\end{theorem}}

\newcommand{\nin}{\noindent}

\newcommand{\co}[1]{} 
\newcommand{\ilfloor}{\downharpoonright\!}
\newcommand{\irfloor}{\!\downharpoonleft}
\newcommand{\ilceil}{\upharpoonright\!}
\newcommand{\irceil}{\!\upharpoonleft}
\newcommand{\ifloor}[1]{\ilfloor #1 \irfloor}
\newcommand{\iceil}[1]{\ilceil #1 \irceil}
\newcommand{\ckp}{{\rm CKP}}
\newcommand{\norm}[1]{\parallel \! #1 \! \parallel}

\newcommand{\argmin}{\operatorname{argmin}}

\newcommand{\p}{\phantom{-}}

\theoremstyle{plain}
\newtheorem{theorem}{Theorem}[section]
\newtheorem{lemma}[theorem]{Lemma}
\newtheorem{definition}[theorem]{Definition}

\newtheorem{example}[theorem]{Example}

\newcommand{\avec}{\mathbf{a}}
\newcommand{\bvec}{\mathbf{b}}
\newcommand{\pvec}{\mathbf{p}}
\newcommand{\rvec}{\mathbf{r}}
\newcommand{\svec}{\mathbf{s}}
\newcommand{\xvec}{\mathbf{x}}
\newcommand{\yvec}{\mathbf{y}}
\newcommand{\cvec}{\mathbf{c}}
\newcommand{\evec}{\mathbf{e}}
\newcommand{\uvec}{\mathbf{u}}
\newcommand{\vvec}{\mathbf{v}}
\newcommand{\zeros}{\mathbf{0}}
\newcommand{\ones}{\mathbf{1}}

\newcommand{\nz}{\boxtimes}
\newcommand{\Lt}{ {\mathcal L} }
\newcommand{\NLt}{ {\mathcal N} }

\begin{document}

\title{\bf Thinner is not Always Better: Cascade Knapsack Problems}
\author{\begin{tabular}{c}
\large Bala Krishnamoorthy\\
\normalsize Washington State University, Vancouver, WA, USA\\
\normalsize \href{mailto:bkrishna@math.wsu.edu}{bkrishna@math.wsu.edu}
\end{tabular}
}
\date{}

\maketitle

\begin{abstract}

  \nin In the context of \bnb (\bbns) for integer programming (IP)
  problems, a direction along which the polyhedron of the IP has
  minimum width is termed a {\em thin} direction. We demonstrate that
  a thin direction need not always be a good direction to branch on
  for solving the problem efficiently. Further, the {\em integer
    width}, which is the number of \bb nodes created when branching on
  the direction, may also not be an accurate indicator of good
  branching directions.

\end{abstract}


\section{Introduction} \label{sec_intro}

Deciding which variables, or combinations of variables, to branch on
is perhaps the most critical step that determines the performance of
any \bnb (\bbns) algorithm to solve integer programs (IPs); see
\cite{AcKoMa2005}, for instance. Most state of the art IP and mixed
integer programming (MIP) solvers use branching on individual
variables (rather than their combinations). Variable selection done by
the default \bb branching on individual variables can be viewed as a
special case of selecting hyperplanes in the more general \bb that
branches on linear combinations of variables.  Branching on
hyperplanes is one of the key steps in Lenstra's seminal algorithm for
solving integer programs (IPs) in fixed dimensions \cite{L83}. This
theoretical algorithm finds a sequence of ``thin'' directions, i.e.,
directions along which the width of the polyhedron of the IP is
``small''. Branching on these thin directions in sequential order will
solve the IP relatively quickly. The thin directions are identified
using lattice basis reduction.

When individual binary variables are present, the default choice to
branch on them (as opposed to finding good combinations of the
variables) is in line with the idea of selecting thin
directions. Recall that a binary variable $x_j$ satisfies $0 \leq x_j
\leq 1$, and hence the width of the polyhedron along the direction of
$x_j$ is not larger than $1$. Motivated in part by the use of thin
directions in Lenstra's algorithm, the author and Pataki previously
created and analyzed \cite{KrPa2009} a general class of
inequality-constrained knapsack feasibility problems called
decomposable knapsack problems (DKPs), whose coefficients have the
form $\avec=\pvec M+\rvec$ for a suitably large positive number
$M$. We showed that DKPs are difficult for ordinary branch-and-bound,
i.e., \bb branching on individual variables (an independent analysis
of the difficulty of \bb on integer knapsacks was presented separately
\cite{K07}).  At the same time, $\pvec$ indicates a thin direction for
the DKPs, and hence the problems are easy when one branches on the
{\em backbone} hyperplane defined by $\pvec \xvec$. The DKPs subsume
several known families of hard IPs, including those proposed by
Jeroslow \cite{J74}, Todd as well as Avis (as attributed to by
Chv\'atal) \cite{C80}, and by Aardal and Lenstra \cite{AL04}. All
these instances suggest that branching on thin directions---as
represented by individual binary variables in typical instances, or by
the backbone constraint given by $\pvec \xvec$ in the case of
DKPs---is typically a good choice for \bb algorithms.

Even when such thin directions are present in an IP instance,
identifying them might not be straightforward. We had previously
proposed \cite{KrPa2009} a simple preconditioning method termed column
basis reduction (CBR) that provides reformulations of general IPs. For
the DKPs, we proved that the thin direction $\pvec\xvec$ is identified
as the last variable in the preconditioned IP given by CBR. Branching
on the last variable solves the problem in one step.

\paragraph{Our Contribution:} We demonstrate that branching on thin
directions need not {\em always} lead to \bb algorithms running
efficiently. On the other hand, branching on certain directions that
are {\em not} thin might lead to quick convergence of \bb algorithms
on certain instances. To this end, we create a generalization of the
DKPs called {\em cascade knapsack problems} (CKPs), which are
inequality-constrained binary knapsack feasibility problems whose
coefficients have the form $\avec = \pvec_1 M_1 + \pvec_2 M_2 +
\rvec$. We show that the width of CKP polytope along $\pvec_1$ is
larger than $1$, while the width along unit directions is equal to
$1$. As such, $\pvec_1$ is not a thin direction. Similarly, $\pvec_2$
is also not a thin direction. Nonetheless, branching on the hyperplane
defined by $\pvec_1 \xvec$ followed by branching on the hyperplane
defined by $\pvec_2 \xvec$ solves the problem quickly, while the
original CKP is hard for ordinary \bbns. A similar behavior is
observed even when we consider {\em integer width}, which is the
number of nodes created when branching on a hyperplane, in place of
width. We also show that this behavior extends to higher order CKPs,
e.g., with $\avec = \pvec_1 M_1 + \pvec_2 M_2 + \pvec_3 M_3 + \rvec$.

We demonstrate that column basis reduction (CBR) is effective in
solving the CKPs quickly in practice. Extending the previous analysis
for DKPs \cite{KrPa2009}, we argue that branching on the collection of
good directions defined by $\pvec_1 \xvec, \pvec_2 \xvec, \pvec_3
\xvec, \dots$ is captured by branching on the last few variables in
the preconditioned IP give by CBR.

\subsection{Related Work \label{ssec_relwork}}

Following Lenstra's algorithm \cite{L83}, Kannan \cite{K83}, and
Lov\'asz and Scarf \cite{LS92} developed similar algorithms for
solving IPs in fixed dimension. Cook et al.~\cite{CRSS93} reported a
practical implementation of the algorithm of Lov\'asz and Scarf
\cite{LS92}, and obtained reductions in the number of nodes in the \bb
tree for solving certain network design problems. At the same time,
finding the branching direction(s) at each node was quite expensive.

Mahajan and Ralphs \cite{MR2010} showed that it is {\cal NP}-hard to
find branching directions that are optimal with respect to width. On
the other hand, Aardal and co-workers studied a basis reduction-based
reformulation technique for equality constrained IPs, which uses the
idea of branching on {\em good} hyperplanes \cite{ABHLS00,AHL00,AL04}.
Specifically, Aardal and Lenstra \cite{AL04} studied a class of
equality-constrained integer knapsack problems whose reformulations
have a specific thin direction, which is also identified by a
reformulation technique. Our preconditioning method CBR
\cite{KrPa2009} applies to more general, i.e., not necessarily
equality-constrained, IPs. The reformulation technique of Aardal et
al.~is subsumed by the more general CBR.

The use of branching on good hyperplanes on more general IPs was
demonstrated by Aardal et al.~\cite{ABHLS00}, who used their basis
reduction-based reformulation technique to solve otherwise
hard-to-solve marketshare problems, which are multiple
equality-constrained binary IPs \cite{CD98,KochEtAl2011}. Louveaux and
Wolsey \cite{LW02} extended these results to a more general class of
IPs with some special structure. In related work, Mehrotra and Li
\cite{ML05} proposed a general framework for identifying branching
hyperplanes for mixed IPs, which also benefits from basis reduction.

Pataki, Tural, and Wong \cite{PaTuWo2010} studied the efficacy of \bnb
on CBR-type reformulations of general IPs. In particular, they proved
an upper bound on the width of the polyhedron of the reformulated IP
along the last unit vector. This bound implies that as the size of
coefficients in the constraint matrix of the original IP increases,
the reformulation of almost all instances generated from a standard
distribution is solved at the root node. This result follows from the
fact that the last unit vector in the reformulation is equivalent to a
direction along which the width of the original polyhedron is small,
and hence branching on the last unit vector solves the problems
easily.

\section{Width, Integer Width, and Branching Directions: Examples} \label{sec_widiwidex}

\bdef \label{def-widiwid}

Given a polyhedron $Q$, and an integer vector $\cvec$, the {\em
  width} and the {\em integer width} of $Q$ in the direction of
$\cvec$ are
\begin{align*}
  \width(\cvec, Q) & = ~\max \{ \cvec\xvec \, | \, \xvec \in Q \} - 
  \min \{ \cvec\xvec \, | \, \xvec \in Q \, \}, \\
  \iwidth(\cvec, Q) & = ~\lfloor \max \{ \cvec\xvec \, | \, \xvec \in
  Q \} \rfloor - \lceil \min \{ \cvec\xvec \, | \, \xvec \in Q \} 
  \rceil + 1.
\end{align*}
Further, an integer direction $\cvec^*$ is termed a {\em thin}
direction of $Q$ if $~\displaystyle \cvec^* \in
~\argmin\limits_{\cvec}\, \width(\cvec,Q)$. The quantity
$\iwidth(\cvec, Q)\,$ is the number of nodes generated by \bnb when
branching on the constraint given by $\cvec\xvec$.

\Edef

We point out that integer width is {\em not} given by requiring the
optima in the definition of width be attained by integer
vectors. Consider the simple 2D example where $Q = \{ \xvec \in \R^2
\, | \, 0 \leq x_1, x_2 \leq 1\}$. For $\cvec = \begin{bmatrix} 1 &
  0 \end{bmatrix}$, we get that $\width(\cvec, Q) = 1$, with the
maximum and minimum being attained by $\xvec = \begin{bmatrix} 1 &
  0 \end{bmatrix}$ and the zero vector, respectively. But
$\iwidth(\cvec, Q) = 2$, since we would consider two branches when
branching on the direction defined by $x_1$, corresponding to $x_1=1$
and $x_1=0$.

We present several instances of integer programming which demonstrate
that thin directions might not always be the best choices for
branching. On the contrary, certain specific directions along which
the width, or even the integer width, of the polytope is {\em larger}
than the minimum (integer) width could help solve the problem quickly
using \bbns.

We consider classes of inequality-constrained binary knapsack
feasibility problems of the form
\[ \{ \xvec \in \zad{n} ~ | ~ \beta' \leq \avec \xvec \leq \beta, 
~\zeros \leq \xvec \leq \ones \}, \]
that have {\em no} integer feasible solutions, and our goal is to {\em
  prove} its integer infeasibility using \bbns. For an integer program
labeled ${\rm (IP)}$ and an integer vector $\cvec$, we denote by
$\width(\cvec, {\rm IP})$ the width of the LP-relaxation of ${\rm
  (IP)}$ in the direction $\cvec$, and similarly denote
$\iwidth(\cvec, {\rm IP})$.

Our use of the term knapsack problem is a generalization of how it is
referred to in most literature (see, e.g., \cite[Section 16.6]{S86}),
where the single constraint has the form $\avec \xvec \leq \beta$ (or
$\avec \xvec = \beta$ in the equality version). It happens to be the
case that we have $\beta' = \beta$ in all the following examples
(\ref{kp1},\ref{kp2},\ref{kp3},\ref{kp4}). But our construction (in
Section \ref{sec_tckp}) is more general, allowing $\beta' <
\beta$. Indeed, several of the larger instances we present (see Table
\ref{tab_w24ckp_n50}) do have $\beta' < \beta$.

\bex \label{ex1iw13ckp} 
  {\rm Consider the following knapsack feasibility problem with $n=12$
    binary variables:
    %
    \beq \tag{\rm{KP1}} \label{kp1}
    \ba{rcl}
    75 \, x_1 +  ~86 \, x_2 + ~97 \, x_3 + ~105 \, x_4  + ~142 \, x_5  + ~153 \, x_6 ~+ \\
    161 \, x_7 + 172 \, x_8 + 209 \, x_9 + 217 \, x_{10} + 228 \, x_{11} + 239 \, x_{12} ~& = & 1023 \\
    0 \leq \, x_j \, \leq 1,~~x_j \in \, \Z,~~j=1,\dots,12. 
    \ea
    \eeq

    \nin There are no integer feasible solutions, and CPLEX 12.6.3.0
    applying branch-and-cut on the original variables without using
    any objective function takes $484$ \bb nodes to prove integer
    infeasibility of this instance (we mention here that {\em all}
    computations presented in this document were done on an Intel PC
    with $8$ cores and a $2.33$ GHz CPU using CPLEX 12.6.3.0 as the
    MIP solver). The knapsack coefficient vector $\avec =
    (75,~86,\dots,239)$ is a thin direction (trivially), as
    $\width(\avec,\mbox{\ref{kp1}})=0$. It can also be checked that
    $\width(\evec_j,\rm{\ref{kp1}})=1$ for all unit vectors $\evec_j$,
    as the maximum and minimum of $x_j$ over (\ref{kp1}) are $1$ and
    $0$, respectively, for all $j \in \{1,\dots,12\}$.

    \medskip
    The knapsack coefficients have the form $\avec = \pvec_1 M_1 + \pvec_2
    M_2 + \rvec$, where
    \beq \label{eq_p1p2r3ckp}
    \ba{rcccccccccccccc}
    \pvec_1 & = & (\p1, &  1, &  1, & 1, &  2, &  2, &  2, &  2, &  3, &  3, &  3, &  3), & \\
    \pvec_2 & = & (\p1, &  2, &  3, & 4, &  1, &  2, &  3, &  4, &  1, &  2, &  3, &  4), & \mbox{and} \\
    \rvec   & = & (-1, &  0, &  1, & -1, &  0, &  1, & -1, &  0, &  1, & -1, &  0, &  1).
    \ea
    \eeq
    $\width(\pvec_1,\rm{\ref{kp1}}) = 13.99 - 12.13 = 1.86$, which is
    larger than the width along any $\evec_j$. Still, we branch on the
    direction of $\pvec_1$ by adding the constraint $\,\pvec_1 \xvec =
    13$. We now get the maximum and minimum of $\pvec_2 \xvec$ for
    this branch as $16.90$ and $16.10$. Thus, branching on the
    hyperplanes defined by $\pvec_1\xvec$ and $\pvec_2\xvec$ in that
    order proves the integer infeasibility of the instance easily, in
    only two \bb nodes. We also point out that $\width(\pvec_2,
    \rm{\ref{kp1}}) = 22.16 - 9.99 = 12.17$, and hence does not
    define a helpful direction to branch on.  }
  \eex
     
  One could argue that while $\width(\pvec_1,\ref{kp1})$ is indeed
  larger than the width along any individual variable,
  $\iwidth(\pvec_1,\ref{kp1})=1$ is in fact smaller than
  $\iwidth(\evec_j,\ref{kp1})=1-0+1=2$. But the next two examples
  illustrate that small integer widths might also not indicate good
  branching directions.

\bex \label{ex2iw23ckp} 
  {\rm Consider the following knapsack feasibility problem with $n=12$
    binary variables:
    %
    \beq \tag{\rm{KP2}} \label{kp2}
      \ba{rcl}
        71 \, x_1 + ~\,82 \, x_2 + ~93 \, x_3 + ~101 \, x_4 + 134 \, x_5 + 145 \, x_6 ~~~~ + \\
        153 \, x_7+ 164 \, x_8 + 197 \, x_9 + 205 \, x_{10} + 216 \, x_{11} + 227 \, x_{12} ~~& = & 981 \\
        0 \leq \, x_j \, \leq 1,~~x_j \in \, \Z,~~j=1,\dots,12. \\
      \ea
    \eeq
    %
    
    There are no integer feasible solutions, and CPLEX 12.6.3.0
    applying branch-and-cut on the original variables without using
    any objective function takes $462$ \bb nodes to prove integer
    infeasibility of this instance. Similar to (\ref{kp1}), we get
    $\width(\avec,\ref{kp2})=0$ and $\width(\evec_j,\ref{kp2})=1$ for
    all unit vectors $\evec_j$.

    The knapsack coefficient vector $\avec = (71,~82, \dots, 227)$ has
    the same structure as in (\ref{kp1}), using the same vectors
    $\pvec_1,\pvec_2,\rvec$ in Equation (\ref{eq_p1p2r3ckp}), except
    with $M_1 = 62$. We observe that $\width(\pvec_1,$ $\ref{kp2}) =
    14.17 - 12.22 = 1.95$ and $\iwidth(\pvec_1,\ref{kp2}) = 2$, which
    is equal to the integer width along each $\evec_j$. Nonetheless,
    branching on the hyperplane defined by $\pvec_1 \xvec$ followed by
    that defined by $\pvec_2 \xvec$ solves the problem in only three
    \bb nodes, as shown in Figure
    \ref{fig_bbtreekp2ckpn12iw2}. Similar to (\ref{kp1}),
    $\width(\pvec_2,\rm{\ref{kp2}})=22.25 - 10.23 = 12.02$, and hence
    $\pvec_2 \xvec$ does not define a good direction to branch on.

    \begin{figure}[ht!]
      \centering
      \includegraphics[scale=0.9]{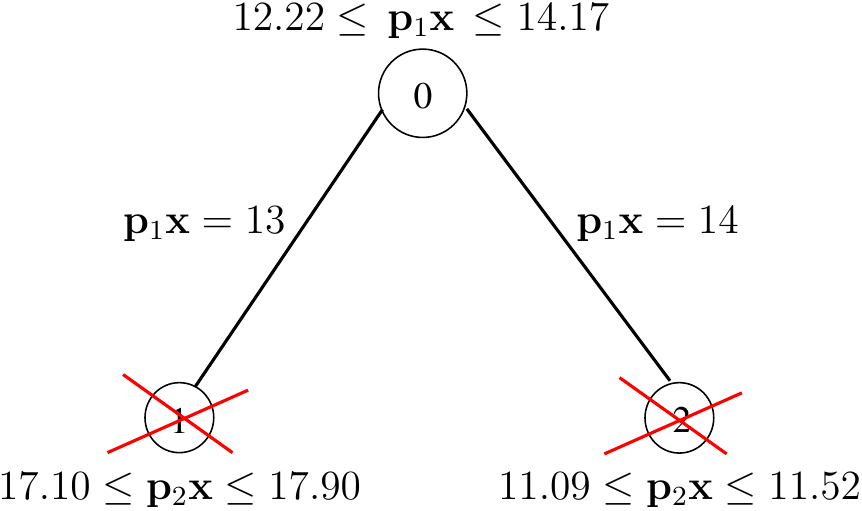}
      \caption{\label{fig_bbtreekp2ckpn12iw2} Branching tree for
        (\ref{kp2}). Crossed out nodes are pruned due to integer
        infeasibility.}
    \end{figure}
  }
  \eex

\bex \label{ex3iw33ckp} 
  {\rm Consider the following knapsack feasibility problem with $n=12$
    binary variables:
    %
    \beq \tag{\rm{KP3}} \label{kp3}
      \ba{rcl}
        57 \, x_1 + 68 \, x_2 + ~79 \, x_3 + ~87 \, x_4 + ~106 \, x_5 + ~117 \, x_6 ~~~~+ \\
        125 \, x_7+ 136 \, x_8 + 155 \, x_9 + 163 \, x_{10} + 174 \, x_{11} + 185 \, x_{12} ~& = & 847 \\
        0 \leq \, x_j \, \leq 1,~~x_j \in \, \Z,~~j=1,\dots,12. \\
      \ea
    \eeq
    %

    There are no integer feasible solutions, and CPLEX 12.6.3.0
    applying branch-and-cut on the original variables without using
    any objective function takes $451$ \bb nodes to prove integer
    infeasibility of this instance. Similar to (\ref{kp1}) and
    (\ref{kp2}), we get $\width(\avec,\ref{kp3})=0$ and
    $\width(\evec_j,$ $\ref{kp3})=1$ for all unit vectors $\evec_j$.

    The knapsack coefficient vector $\avec = (57,~68, \dots, 185)$ has
    the same structure as in (\ref{kp1}), using the same vectors
    $\pvec_1,\pvec_2,\rvec$ in Equation (\ref{eq_p1p2r3ckp}), except
    with $M_1 = 48$.
    We observe that $\width(\pvec_1,$ $\ref{kp3}) = 15.22 - 12.86 =
    2.36$ and $\iwidth(\pvec_1,\ref{kp3}) = 3$, which is strictly
    greater than the integer width along each $\evec_j$. Nonetheless,
    branching on the hyperplane defined by $\pvec_1 \xvec$ followed by
    that defined by $\pvec_2 \xvec$ solves the problem in only four
    \bb nodes, as shown in Figure
    \ref{fig_bbtreekp3ckpn12iw3}. Similar to (\ref{kp1}) and
    (\ref{kp2}), $\width(\pvec_2,\rm{\ref{kp3}})=22.90 - 11.62 =
    11.28$, and $\pvec_2 \xvec$ does not define a good direction to
    branch on.
    \begin{figure}[ht!]
      \centering
      \includegraphics[scale=0.85]{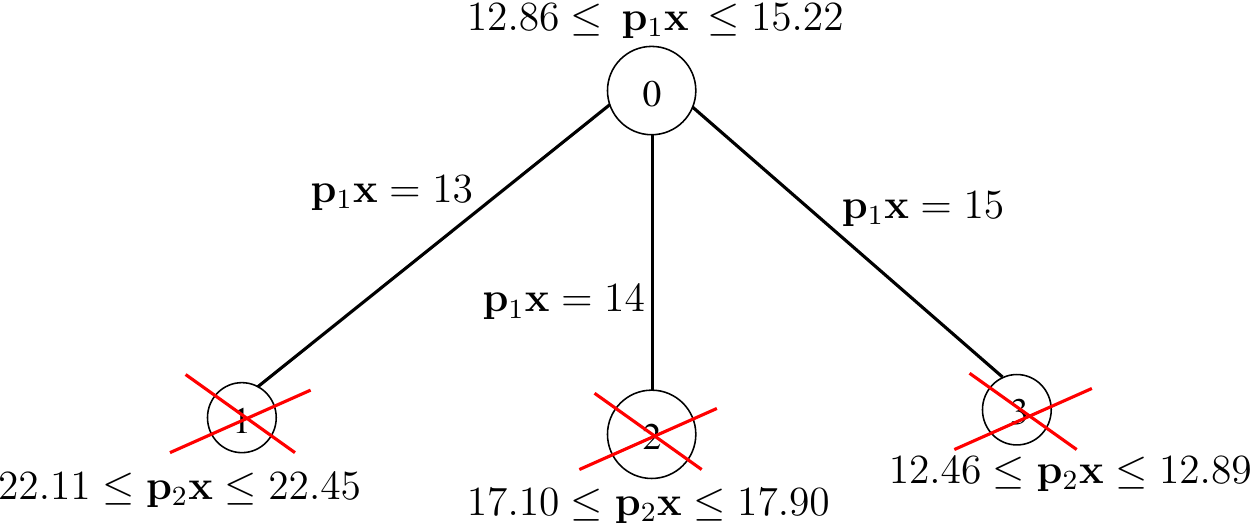}
      \caption{\label{fig_bbtreekp3ckpn12iw3} Branching tree for
        (\ref{kp3}). Crossed out nodes are pruned due to integer
        infeasibility.}
    \end{figure}

  }
  \eex

  These three examples have a common property: the effect of branching
  on the hyperplane defined by $\pvec_1 \xvec$ {\em cascades} down to
  the subproblems created in this process, which are all pruned by
  branching subsequently on the hyperplane defined by
  $\pvec_2\xvec$. Hence we call them {\em cascade knapsack problems}
  (CKPs). The next example illustrates that the cascading effect could
  be observed over multiple $\pvec_i \xvec$ directions, i.e.,
  branching on the hyperplane defined by $\pvec_1\xvec$ cascades down,
  and then the effect of branching on the hyperplane defined by
  $\pvec_2\xvec$ cascades down to the next levels of nodes, followed
  by branching on the hyperplane defined by $\pvec_3\xvec$, and so on.


%
\bex \label{ex4iw24ckp} 
  {\rm Consider the following knapsack feasibility problem with $n=12$
    binary variables:
    \beq \tag{\rm{KP4}} \label{kp4}
      \ba{rcl}
        723 \, x_1 + ~799 \, x_2 + ~875 \, x_3 ~+ ~981 \, x_4 ~+ ~1285 \, x_5 + ~1361 \, x_6 ~~~+ \\
        1467 \, x_7 + 1587 \, x_8 + 1847 \, x_9 + 1953 \, x_{10} + 2029 \, x_{11} + 2116 \, x_{12} ~~& = & 9312 \\
        0 \leq \, x_j \, \leq 1,~~x_j \in \, \Z,~~j=1,\dots,12. \\
      \ea
    \eeq

    The knapsack coefficients have the form $\, \avec = \pvec_1 M_1 +
    \pvec_2 M_2 + \pvec_3 M_3 + \rvec \,$ for the same set of vectors
    $\pvec_1, \pvec_2,$ and $\rvec$ given in Equation
    (\ref{eq_p1p2r3ckp}), with $\pvec_3 = (\,5,~ 3,~ 1,~ 2,~ 4,~ 2,~
    3,~ 5,~ 3,~ 4,~ 2,~ 1\,)$, and $M_1 = 572, M_2 = 97$, and $M_3 =
    11$. The vector $\pvec_3$ is chosen to be linearly independent of
    $\pvec_1$ and $\pvec_2$ here. This knapsack problem is also
    integer infeasible, and CPLEX 12.6.3.0 branching on the original
    variables takes $491$ \bb nodes to prove integer
    infeasibility. Further, $\width(\avec,\ref{kp4})=0$ and
    $\width(\evec_j,\ref{kp4})=1$ for all $j$, as in the previous
    instances. At the same time, we can solve this problem easily if
    we branch on the hyperplanes defined by $\pvec_1\xvec,
    \pvec_2\xvec,$ and $\pvec_3\xvec$, in that order. The branching
    tree is shown in Figure \ref{fig_bbtrex4ckpn12iw2}, and the
    problem is solved in seven \bb nodes.

    \begin{figure}[ht!]
      \centering
      \vspace*{0.1in}
      \includegraphics[scale=0.84]{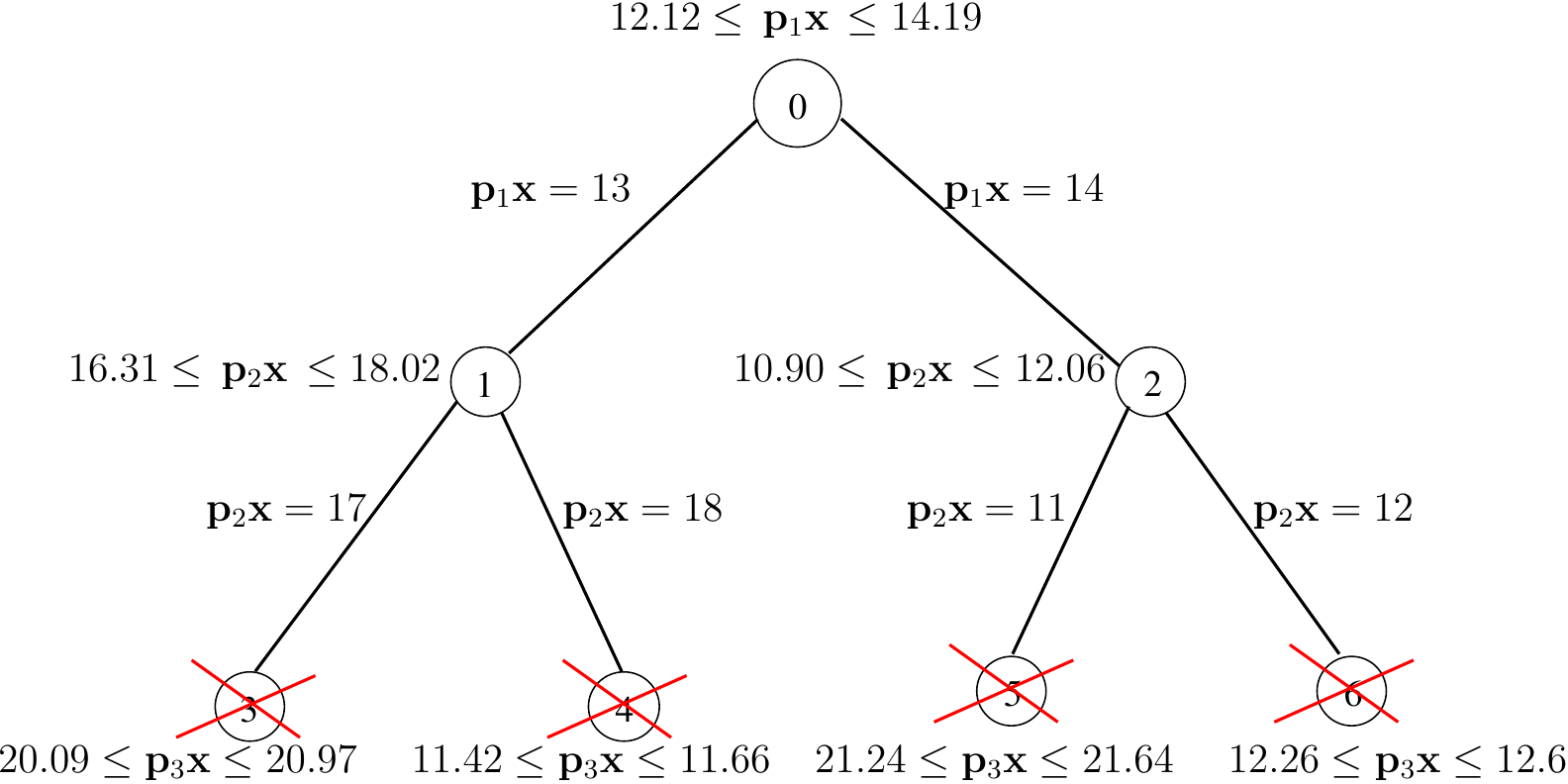}
      \caption{\label{fig_bbtrex4ckpn12iw2} Branching tree for
        (\ref{kp4}). Crossed out nodes are pruned due to integer
        infeasibility.}
    \end{figure}
  }
  \eex

We could present a similar instance with $\, \avec = \pvec_1 M_1 +
\pvec_2 M_2 + \pvec_3 M_3 + \rvec \,$ where the integer width of its
LP relaxation along $\pvec_1$ is $3$, similar to (\ref{kp3}). Just as
in (\ref{kp4}), we solve this problem quickly if we branch on the
hyperplanes defined by $\pvec_1\xvec, \pvec_2\xvec,$ and
$\pvec_3\xvec$, in that order.

\section{Procedure for creating a CKP \label{sec_tckp}}

We present a procedure (Figure \ref{fig_procckp}) to create CKP
instances with the structure illustrated in Example \ref{ex1iw13ckp}
(\ref{kp1}), where $\avec = \pvec_1 M_1 + \pvec_2 M_2 + \rvec$ and
branching on the hyperplane defined by $\pvec_1\xvec$ followed by
branching on the hyperplane defined by $\pvec_2\xvec$ solves the
problem easily. The remaining examples (\ref{kp2}, \ref{kp3},
\ref{kp4}) as well as the instances reported in Section
\ref{sec_comp4ckps} were generated by appropriate modifications of
this procedure. We use the following notation in the description.

\bdef \label{def-ifloorceil}
For $x \in \R$, we define
\[ \ifloor{x} = \left\{ \ba{c@{\mbox{ if }}l}
\lfloor x \rfloor & x \not\in \Z, \\
x - 1 & x \in \Z; \ea \right.~~~\mbox{and}~~~
\iceil{x} = \left\{  \ba{c@{\mbox{ if }}l}
\lceil  x \rceil  & x \not\in \Z, \\
x + 1 & x \in \Z. \ea \right.
\]
\Edef

\paragraph[other]{Other notation:} We denote the vector of 
ones by $\ones$, and the box with upper bound $\ones$ by $B_{\ones} =
\{ \, \xvec \in \R^n \, | \, \zeros \leq \xvec \leq \ones \}$.
We denote by $\zad{n}_+$ the set of all $n$-vectors with positive integer entries.

\medskip

\begin{figure}[ht!] 
\label{proc1lab}
\framebox[6.5in]{\parbox{6.3in}{
{\sc Procedure CKP}
\begin{center}
\begin{tabular}{ll}
{\tt Input:} & Vectors $\pvec_1, \pvec_2, \rvec$ with $\pvec_1,
\pvec_2 \in \zad{n}_+$, $\rvec \in \zad{n}$, and $\pvec_1,
\pvec_2$ linearly independent.  \\
%
{\tt Output:} & $M_1, M_2, \beta', \beta$ forming the CKP instance
with $\avec = \pvec_1 M_1 + \pvec_2 M_2 + \rvec$.
\end{tabular}
\end{center}

\begin{enumerate}
\item {\tt Choice of $k_1, k_2$:} 
  \beq \label{eq_maxminpix}
\hspace*{-0.2in}
\ba{c}
\mbox{Find}~\gamma_1 = \max\{\pvec_1 \xvec \,|\,\xvec \in B_{\ones}\},~ 
\delta_1 = \min\{\pvec_1 \xvec \,|\,\xvec \in B_{\ones}\}; ~\mbox{set}~ 
k_1 = \lfloor\, (\ifloor{\gamma_1} + \iceil{\delta_1})/2 \, \rceil.   \\
\vspace*{-0.1in} \\
\mbox{Find } \gamma_2 = \max\{\pvec_2 \xvec ~|~\pvec_1 \xvec = k_1 + 1, ~\xvec \in B_{\ones}\},~\\
  \mbox{ and } ~\,\delta_2 = \min\,\{\pvec_2 \xvec ~|~\pvec_1 \xvec = k_1 + 1, ~\xvec \in B_{\ones}\};~~ \\
 ~~~~~~~~\mbox{set} ~  k_2 = \lfloor\, (\ifloor{\gamma_2} + \iceil{\delta_2})/2 \, \rceil.  
\ea
\eeq
\vspace*{-0.05in}
\item {\tt Choice of $M_2, \beta'_2, \beta_2$:} 
\beq \label{eq_eta2mu2}
\ba{rcl}
\mbox{Find }~\eta_2 & = & \max\{\rvec \xvec ~|~\pvec_1 \xvec = k_1 + 1, \pvec_2 \xvec \leq k_2,~\xvec \in B_{\ones}\}, \\
~\mbox{ and }\mu_2 & = & \min\,\{\rvec\xvec~|~\pvec_1 \xvec = k_1 + 1, \pvec_2 \xvec \geq k_2+1,~\xvec \in B_{\ones} \}.
\ea
\eeq

\vspace*{-0.1in}
\beq \label{eq_M2def}
\mbox{If}~
\iceil{\eta_2} - \ifloor{\mu_2} \geq 1,~\mbox{set}~
M_2 = \iceil{\eta_2} - \ifloor{\mu_2}, ~\mbox{and set}~~
\avec_2 = \pvec_2 M_2 + \rvec.
\eeq

\vspace*{-0.15in}
\beq \label{eq_eta2prmu2pr}
\ba{rcl}
\mbox{Find }~\eta'_2 & = & \max\,\{\avec_2 \xvec~|~\pvec_1 \xvec =
 k_1+1,\, \pvec_2 \xvec \leq k_2,~\xvec \in B_{\ones}\}, \\
\vspace*{-0.15in} \\
\mbox{ and }\mu'_2 & = & \min\,\,\{\avec_2 \xvec~|~\pvec_1 \xvec = k_1+1,
 ~\pvec_2 \xvec \geq k_2+1,~\xvec \in B_{\ones} \}.
\ea
\eeq

\vspace*{-0.1in}
\beq \label{eq_beta2s}
\mbox{If } \iceil{\eta'_2} \leq \ifloor{\mu'_2}~~\mbox{set}~
\beta'_2 = \iceil{\eta'_2}~\mbox{ and }~\beta_2 = \ifloor{\mu'_2}.
\eeq

\medskip
\item {\tt Choice of $M_1, \beta'_1, \beta_1,$: \ }
  \beq \label{eq_eta1mu1}
    \ba{c}
    \mbox{Find }~\eta_1 ~=~ \max\,\{\avec_2\xvec ~|~ \pvec_1 \xvec
    \leq k_1,~\xvec \in B_{\ones}\},~~~~~~\mbox{ and }~\, \\
    \mu_1 ~ = ~ \min\,\,\{\avec_2\xvec ~|~ \pvec_1 \xvec \geq
    k_1+2,~\xvec \in B_{\ones} \}.
    \ea
  \eeq

\vspace*{-0.2in}
\beq \label{eq_M1}
\mbox{Set }~ M_1 = \iceil{\max\,(\,\eta_1-\beta'_2, \, \beta_2 - 
\mu_1 ,\, M_2)}.
\eeq

\vspace*{-0.15in}
\beq \label{eq_beta1s}
\mbox{Set } ~\avec_1 = \pvec_1 M_1 + \avec_2,~ 
~\beta'_1 = \beta'_2 + (k_1+1) M_1,~\mbox{ and }
~\beta_1  = \beta_2  + (k_1+1) M_1. \mbox{\hspace*{0.2in}}
\eeq

\medskip
\item {\tt Output instance:}  \ Set $\beta'=\beta'_1,\,\beta = \beta_1, \, 
\avec = \avec_1$, and return $(M_1,M_2, \beta', \beta, \avec)$.
\end{enumerate}
}}
\caption{\label{fig_procckp} Procedure to generate CKPs. }
\end{figure}

\vspace*{0.1in}
We now present two lemmas describing the Procedure (Figure
\ref{fig_procckp}). Notice that we refer to both the final knapsack
problem and its LP relaxation as (CKP), with the exact choice evident
from the context.

\ble
\label{lem-ckpp1k11nonempty}

If the quantities defined in Equations (\ref{eq_M2def}),
(\ref{eq_beta2s}), and (\ref{eq_M1}) exist, then \rm{(CKP)}$\wedge \,
(\pvec_1 \xvec = k_1+1) \, \neq \, \emptyset$.  \ele
\begin{proof}
Let $\xvec_1$ and $\xvec_2$ attain the optima for the LPs in Equation
(\ref{eq_eta1mu1}) (the maximum and the minimum, respectively), based
on which $M_1$ is defined in Equation (\ref{eq_M1}). We get
\beqast
\max\{\avec_1\xvec ~|~\pvec_1\xvec \leq k_1, \xvec \in B_{\ones}\} & \leq & 
k_1 M_1 + \max\{\avec_2\xvec~|~\pvec_1\xvec \leq k_1, \xvec \in B_{\ones}\} 
= k_1 M_1 + \eta_1\\
 & = & [ (k_1+1)M_1 + \beta'_2 ] - \, (M_1 -(\eta_1 - \beta'_2)) \\
 & = & \beta'_1 - \, (M_1 -(\eta_1 - \beta'_2)) \\
 & < & \beta'_1.
\eeqast
The third equality follows from the definition of $\beta'_1$ in
Equation (\ref{eq_beta1s}), and the last strict inequality follows
from the definition of $M_1$ in Equation (\ref{eq_M1}). Similarly, we
get
\beqast
\min\{\avec_1\xvec ~|~\pvec_1\xvec \geq k_1+2, \xvec \in B_{\ones}\} & \geq & 
(k_1+2) M_1 + \min\{\avec_2\xvec~|~\pvec_1\xvec \geq k_1+2,
\xvec \in B_{\ones}\} \\
 & = & (k_1+2) M_1 + \mu_1\\
 & = & [ (k_1+1)M_1 + \beta_2 ] + \, (M_1 -(\beta_2 - \mu_1)) \\
 & = & \beta_1 + \, (M_1 -(\beta_2 - \mu_1)) \\
 & > & \beta_1.
\eeqast
Hence $\xvec' \in \{ \xvec ~|~\beta'_1 \leq \avec_1 \xvec \leq
\beta_1, \, \xvec \in B_{\ones} \}\,$ implies that $\,k_1 < \pvec_1
\xvec' < k_1+2$. Further, since $\avec_2$ and $\pvec_1$ have positive
entries, the structure of the LPs in Equation (\ref{eq_eta1mu1})
implies that we will have $\pvec_1 \xvec_1 = k_1$ and $\pvec_1 \xvec_2
= k_1 + 2$. Thus we can find a constant $\lambda \in (0,1)$ such that
$\xvec' = \lambda \xvec_1 + (1-\lambda) \xvec_2\,$ satisfies
$\,\beta'_1 \leq \avec_1 \xvec' \leq \beta_1$ with $\pvec_1 \xvec' =
k_1 + 1$. This result implies when we branch on the hyperplane defined
by $\pvec_1 \xvec$ in the original (CKP) problem, the branch created
by setting $\pvec_1 \xvec = k_1 +1$ cannot be pruned due to LP
infeasibility.
\end{proof}

\ble
\label{lem-p1p2provinfeas}
If the quantities defined in Equations (\ref{eq_M2def}),
(\ref{eq_beta2s}), and (\ref{eq_M1}) exist, then the integer
infeasibility of \rm{(CKP)} is proven by $~(\pvec_1 \xvec = k_1+1) \,
\bigwedge \, (\pvec_2 \xvec \leq k_2 \,\vee \, \pvec_2 \xvec \geq
k_2+1)$.  \ele
\begin{proof}
Let $\xvec_3$ and $\xvec_4$ attain the optima for the LPs in Equation
(\ref{eq_eta2prmu2pr}) (the maximum and the minimum, respectively).
Since $\,\eta'_2 < \beta'_2 \leq \beta_2 < \mu'_2,\,$ and due to the
way $\avec_1, \beta'_1,$ and $\beta_1$ are defined (in Equation
(\ref{eq_beta1s})), there is some $\lambda' \in (0,1)$ such that
$\xvec'' = \lambda' \xvec_3 + (1-\lambda') \xvec_4$ and $\xvec'' \in
\{ \xvec ~|~ \beta'_1 \leq \avec_1 \xvec \leq \beta_1,\, \pvec_1 \xvec
= k_1+1,\, \xvec \in B_{\ones} \}$. Since $\avec_2$ and $\pvec_2$ both
have positive entries, the structure of LPs in Equation
(\ref{eq_eta2prmu2pr}) implies we have $\pvec_2 \xvec_3 = k_2$ and
$\pvec_2 \xvec_4 = k_2+1$. Further,
\[
 \begin{array}{ll}
  & \xvec'' \in \{ \xvec \,|\, \beta'_1 \leq \avec_1 \xvec \leq
  \beta_1,\, \pvec_1 \xvec = k_1+1,\, \xvec \in B_{\ones} \} \\
  ~\Rightarrow  & \xvec'' \in \{ \xvec \,|\, \beta'_2 \leq \avec_2
  \xvec \leq \beta_2,\, \pvec_1 \xvec = k_1+1,\, \xvec \in B_{\ones}
  \}  \\
  ~\Rightarrow & k_2 < \pvec_2 \xvec'' < k_2+1,
\end{array}
\]
as $\xvec'' = \lambda' \xvec_3 + (1-\lambda') \xvec_4$ must hold for
some $\lambda' \in (0,1)$.  Thus (CKP)$\wedge \, (\pvec_1 \xvec =
k_1+1) \, \neq \, \emptyset$ as already shown by Lemma
\ref{lem-ckpp1k11nonempty}, and (CKP)$\wedge~(\pvec_1 \xvec = k_1+1)
\, \bigwedge \, (\pvec_2 \xvec \leq k_2 \,\vee \, \pvec_2 \xvec \geq
k_2+1)\, = \, \emptyset$.
%
\end{proof}

\medskip
The preceding two lemmas show that branching on the hyperplane defined
by $\pvec_1 \xvec$ followed by that defined by $\pvec_2\xvec$ proves
the integer infeasibility of (CKP). While Lemma
\ref{lem-ckpp1k11nonempty} shows the existence of some $\xvec' \in$
(CKP) satisfying $\pvec_1 \xvec' = k_1 + 1$, we have not explicitly
shown that $\width(\pvec_1,\ckp) > 1$. In fact, the Procedure (in
Figure \ref{fig_procckp}) is not guaranteed to work on {\em every}
choice of the input vectors. For instance, we might not get a valid
choice for $M_2$ in Equation (\ref{eq_M2def}), or for $\beta'$ and
$\beta$ in Equation (\ref{eq_beta2s}). Further, while the Procedure
assumes only that $\pvec_1$ and $\pvec_2$ are linearly independent, we
might not get the required structure for the CKP if they are ``too
close to each other'', e.g., when $n$ is large and the two vectors
differ in only one or two entries.

To create instances with $\iwidth(\pvec_1,\ckp)=2$ as in (\ref{kp2}),
we could replace the constraint $\pvec_1 \xvec = k_1+1$ in Equations
(\ref{eq_maxminpix}), (\ref{eq_eta2mu2}), and (\ref{eq_eta2prmu2pr})
with $~k_1+1 \leq \pvec_1 \xvec \leq k_1+2$, and the constraint
$\pvec_1 \xvec \geq k_1+2$ in Equation (\ref{eq_eta1mu1}) with
$\pvec_1 \xvec \geq k_1+3$. Alternatively, we could make only the
latter change (to $\pvec_1 \xvec \geq k_1+3$), while sticking with
$\pvec_1 \xvec = k_1+1$ in the former Equations, as we did in
computational tests reported in the following Section.

\vspace*{0.1in}
\nin {\bf Example \ref{ex1iw13ckp} continued:} We created the instance
(\ref{kp1}) using the procedure described in Figure
\ref{fig_procckp}. We obtained $\gamma_1 = 24$ and $\delta_1 = 0$,
setting $k_1=12$, and then $\gamma_2 = 23,\, \delta_2 = 9$, setting
$k_2=16$. The remaining parameter values are obtained as
follows. $\eta_2 = 4, \, \mu_2 = -4,$ giving $M_2 = 10$, and $\eta'_2
= 164, \, \mu'_2 = 166$, and hence $\beta'_2 = \beta_2 = 165$. We
finally get $\eta_1 = 221, \, \mu_1 = 100$, giving $M_1=66$, and
$\beta'_1 = \beta_1 = 1023$.

\section{Computational Tests on CKPs} \label{sec_comp4ckps}

We used a modification of {\sc Procedure CKP} to generate ten
instances of CKP with the structure illustrated by Example
\ref{ex4iw24ckp} (\ref{kp4}), i.e., with $\avec_1 = \pvec_1M_1 +
\pvec_2M_2 + \pvec_3M_3+ \rvec$, for $n=50$. To keep the knapsack
coefficients relatively small, the entries of $\pvec_1,\pvec_2,$ and
$\pvec_3$ are chosen randomly from $\{1,2,3\}$ and those of $\rvec$
from $\{-1,0,1\}$, such that no two $a_j$'s are identical. We tried to
solve the original formulations of the CKP, the CKP with
$\pvec_1\xvec$ fixed, and the CKP with $\pvec_1\xvec$ and
$\pvec_2\xvec$ fixed. All calculations are done on an Intel PC with
$8$ cores and a $2.33$ GHz CPU. As an MIP solver, we used CPLEX
12.6.3.0. For feasibility versions of integer programs, the sum of the
variables is used as the dummy objective function.  Ideally, these
problems are expected to become easier when $\pvec_1\xvec$, and then
$\pvec_2\xvec$, are fixed. At the same time, all the subproblems
obtained by fixing $\pvec_1\xvec = k_1+1$ and then $\pvec_2\xvec =
k_2+1$ still remain relatively hard---they all remain unsolved after
one hour of computational time.  For the record, the number of \bb
nodes examined within this time for all the runs was 41 $\pm$ 12
million (mean $\pm$ std.~dev.). The details of the computations are
provided in Table \ref{tab_w24ckp_n50}.  Notice that the CKP instances
remain relatively hard for ordinary \bb even after branching on both
the hyperplanes defined by $\pvec_1\xvec$ and $\pvec_2\xvec$.

\begin{table}[ht!]
  \centering
    \parbox{6.5in}{\caption{ \label{tab_w24ckp_n50} Statistics for
        CKP instances of size $n=50$ with $\iwidth(\pvec_1|\ckp)=2$,
        $\iwidth(\pvec_2|(\ckp) \wedge \pvec_1\xvec=k_1+1)=2$ and
        $\iwidth(\pvec_3|(\ckp) \wedge \pvec_1\xvec=k_1+1 \wedge
        \pvec_2\xvec=k_2+1)=0$. $a_{\min}$ and $a_{\max}$ give the
        smallest and largest knapsack coefficients.  $w_i$ gives
        $\width(\pvec_i,\ckp)$ for $i=1,2,3$. $\w_{j1}$ gives
        $\width(\pvec_j, (\ckp)\wedge \pvec_1 \xvec = k_1+1)$ for $j =
        2, 3$. $\w_{312}$ gives $\width(\, \pvec_3, (\ckp) \wedge
        \pvec_1 \xvec = k_1+1 \,\wedge \, \pvec_2 \xvec = k_2+1)$. For
        each instance, the original CKP, the CKP after fixing
        $\pvec_1\xvec=k_1+1$, as well as the CKP after fixing both
        $\pvec_1\xvec=k_1+1$ and $\pvec_2\xvec=k_2+1$ were all
        unsolved after the one hour time limit. CBR gives the number
        of \bb nodes examined to solve the CBR-based reformulation,
        which was integer infeasible for every instance. The times
        taken to obtain the CBR-reformulation and to solve it were
        negligible -- each step took less than $1$ second. }}

    \vspace*{0.1in}
    \begin{tabular}{|r|c|c|c|c|c|c|c|c|c|c|c|} \hline 
      \ \ \  & \multicolumn{4}{|c|}{CKP numbers} &  \multicolumn{6}{|c|}{CKP widths}  & {CBR} \\ \hline 
      \# &  $a_{\min}$ &  $a_{\max}$ &  $\beta'$ &  $\beta$  &  $\w_{1}$ &  $\w_{2}$ &  $\w_{3}$ &  $\w_{21}$ &  $\w_{31}$ &  $\w_{312}$ &   BB \\ \hline 
  1 & 13354 & 26674 & 424920 & 424921 & 2.005 & 42.02 & 41.89 & 2.443 & 38.92 & 0.946 &  7 \\ \hline
  2 & 12251 & 24467 & 367732 & 367733 & 2.073 & 40.90 & 43.07 & 2.421 & 38.84 & 0.944 & 11 \\ \hline
  3 & 14456 & 28877 & 416513 & 416514 & 2.050 & 42.91 & 40.35 & 2.083 & 37.37 & 0.944 &  8 \\ \hline
  4 & 14549 & 28490 & 461490 & 461491 & 2.033 & 42.95 & 42.17 & 2.037 & 38.57 & 0.941 &  2 \\ \hline
  5 & 15375 & 30716 & 457004 & 457005 & 2.007 & 42.01 & 41.82 & 2.011 & 39.60 & 0.943 & 10 \\ \hline
  6 & 11234 & 21946 & 326617 & 326618 & 2.077 & 38.93 & 40.86 & 2.449 & 39.02 & 0.943 & 20 \\ \hline
  7 & 11306 & 22578 & 336621 & 336622 & 2.077 & 38.97 & 39.74 & 2.369 & 37.73 & 0.943 & 16 \\ \hline
  8 & 14696 & 29358 & 437657 & 437658 & 2.054 & 44.96 & 41.80 & 2.161 & 38.74 & 0.943 & 10 \\ \hline
  9 & 15722 & 31407 & 453190 & 453191 & 2.050 & 44.91 & 42.20 & 2.141 & 37.95 & 0.947 & 12 \\ \hline
 10 & 15145 & 30255 & 466572 & 466573 & 2.036 & 42.95 & 39.77 & 2.011 & 37.61 & 0.944 & 14 \\ \hline
    \end{tabular} 
  \end{table}

All these instances of (CKP) are integer infeasible. They have an
integer width of $2$ along $\pvec_1$, and for one of the two
subproblems created by branching on the hyperplane defined by
$\pvec_1\xvec$, the integer width along $\pvec_2$ is $2$ as well. The
other subproblem is not guaranteed to have this property. In four of
the instances, the other subproblem is solved by branching on the
hyperplane defined by $\pvec_2 \xvec$ (integer width along $\pvec_2$
is zero), while for the remaining six instances, the integer width
along $\pvec_2$ is $1$. For all ten instances, the integer
infeasibility of the subproblems created by branching on the
hyperplane defined by $\pvec_2\xvec$ (if any) is proven by branching
on the hyperplane defined by $\pvec_3\xvec$ in the next level. As
mentioned previously, this modified version of {\sc Procedure CKP}
cannot be guaranteed to work for every choice of problem
parameters. Still, we used the modified procedure as a guideline to
search for appropriate parameters that generated instances with the
desired structure. As part of the computational tests, we tried to
solve the original CKP, the original CKP with $\pvec_1\xvec$ fixed,
and also the original CKP with $\pvec_1\xvec$ and $\pvec_2\xvec$
fixed.

\medskip
\nin The instances are available online at
\href{http://www.math.wsu.edu/faculty/bkrishna/CKP/}{http://www.math.wsu.edu/faculty/bkrishna/CKP/}.

\subsection{Column Basis Reduction and CKPs} \label{ssec_cbrckp}

We applied the reformulation technique termed column basis reduction
(CBR) \cite{KrPa2009} on the CKP instances. CBR is a simple
preconditioning method for IP feasibility that replaces the problem
\begin{equation}
  \label{eq-CBR}
  \{\xvec \in \Z^n~|~\bvec' \leq A \, \xvec \leq \bvec \}
   ~~\mbox{ with }~~
  \{\yvec \in \Z^n~|~\bvec' \leq AU \, \yvec \leq \bvec \},
\end{equation}
where $U$ is a unimodular matrix computed using basis reduction (BR)
applied on $A$, which makes the columns of $A$ short, i.e., have small
euclidean norms, and nearly orthogonal. The variables $\yvec$ in the
reformulation and the original variables $\xvec$ are related as $\yvec
= U^{-1} \xvec$. Standard methods of BR include the
Lenstra-Lenstra-Lov\'asz (LLL) reduction \cite{LLL82}, which runs in
time polynomial time, and versions of Korkine-Zolotarev (KZ) reduction
including block-KZ or BKZ reduction \cite{Sh87}, which results in a
higher quality of reduction but runs in polynomial time only when the
dimension is fixed.

For basis reduction calculations on the CKP instances, we used the BKZ
reduction algorithm with the number of columns in the matrix used as
the block-size, using the subroutines from the Number Theory Library
(NTL) version 9.4.0 \cite{NTL} with GMP version 6.0.0. The CBR
reformulations of each instance reported in Table \ref{tab_w24ckp_n50}
is solved in a few \bb nodes in less than one second of computational
time.

\bigskip
We had previously analyzed the efficacy of CBR on DKPs
\cite{KrPa2009}, where the quantities in the problems specified in the
Inequality (\ref{eq-CBR}) are specified as follows.
\[
A = \begin{bmatrix} \avec \\ I \end{bmatrix} ~=~ \begin{bmatrix}
  \pvec M + \rvec \\ I \end{bmatrix},
   ~\mbox{ where }~I~\mbox{ is the identity matrix,}~
  \bvec' = \begin{bmatrix} \beta' \\ \zeros \end{bmatrix},~
  \mbox{ and }~
  \bvec  = \begin{bmatrix} \beta \\ \uvec \end{bmatrix}~
  \mbox{ for }~
  \uvec \in \zad{n}_+.
\]
Under some assumptions on the size of $M$ compared to the norms of
$\pvec$ and $\rvec$, we showed that branching on the last few
variables in the CBR reformulation is equivalent to branching on the
hyperplane defined by $\pvec \xvec$ in the original DKP. We now
present a more general analysis which suggests a similar behavior for
the CBR reformulation of CKPs, i.e., branching on the important
directions defined by $\pvec_1\xvec, \pvec_2\xvec, \pvec_3\xvec,
\dots$ is captured by branching on the last few individual variables
in the reformulation.

We first introduce some definitions and notation related to BR. Given
a matrix $B \in \R^{m \times n}$ with $m \geq n$ and linearly
independent columns, the {\em lattice} generated by the columns of $B$
is $\, \Lt (B) = \{ B \xvec ~|~ \xvec \in \Z^n \}$, i.e., the set of
all integer combinations of columns of $B$. The $k$-th successive
minimum of $\Lt (B)$ is
\[
  \Lambda_k (\Lt (B)) = \min \{ t \, | \, \exists k \mbox{ linearly
    independent vectors in } \Lt (B) \mbox{ with norm at most } t \}.
\]
Suppose there is a constant $c_n$ that depends only on $n$ with the
following property: if $\tilde{\bvec}_j$ denote the columns of
$\tilde{B} = BU$ computed by BR, then
\[
\max \,\{ \norm{\tilde{\bvec}_1},  \dots, \norm{\tilde{\bvec}_{\ell}} \}
\leq c_n \Lambda_{\ell} (\Lt (B))~\mbox{ for } \ell=1,\dots,n.
\]
Then $c_n$ is termed the {\em strength} of BR; the smaller the value
of $c_n$, the more reduced are the columns of $BU$. LLL reduction has
strength $c_n = 2^{(n-1)/2}$ while KZ reduction has strength $c_n =
\sqrt{n}$ (see \cite{Sh87}, for instance). Finally, the {\em kernel
  lattice} or {\em null lattice} of the columns of $B$ is $\NLt (B) =
\{ \xvec \in \Z^n \, | \, B \xvec = \zeros\}$.

\medskip
We consider a general CKP of the form $\,\{ \xvec \in \zad{n} ~ | ~
\beta' \leq \avec \xvec \leq \beta, ~\zeros \leq \xvec \leq \ones
\}\,$ where the coefficient vector has the structure $\avec = \pvec_1
M_1 + \pvec_2 M_2 + \dots + \pvec_t M_t + \rvec$ for positive linearly
independent vectors $\pvec_i \in \zad{n}_+, i=1,\dots,t\,$ and $\rvec
\in \zad{n}$ (we assume $t \geq 2$). The multipliers $M_i$ satisfy
$M_1 > M_2 > \dots > M_t$.  The CBR reformulation of this CKP is of
the form $\,\{\yvec \in \zad{n} ~|~ \beta' \leq (\avec U) \yvec \leq
\beta, ~\zeros \leq U \yvec \leq \ones \}$, where $U$ is the
unimodular matrix obtained by applying BR on the matrix
\begin{equation}
  \label{eq-Ackp}
  A = \begin{bmatrix} \avec \\ I \end{bmatrix} ~ = ~
  \begin{bmatrix}
    \pvec_1 M_1 + \pvec_2 M_2 + \dots + \pvec_t M_t + \rvec \\
    I \end{bmatrix}.
\end{equation}
We denote by $\tilde{A} = AU$ the matrix resulting from applying BR on
$A$ in Equation (\ref{eq-Ackp}). With $P \in \zad{t \times n}_+$
denoting the matrix obtained by stacking the rows $\pvec_1, \dots,
\pvec_t$ vertically in that order, we similarly denote $\tilde{P} =
PU$, $\tilde{\avec} = \avec U$, as well as $\tilde{\pvec}_i = \pvec_i
U$ for $i=1,\dots,t$.  For $1 \leq r \leq s \leq n$, we denote the
subset of $r$-th to $s$-th entries of the $i$-th row of $\tilde{P}$
(equivalently of $\tilde{\pvec}_i$) by $\tilde{P}_{i,r:s}$ (or
$\tilde{\pvec}_{i,r:s}$). The following theorem describes why CBR
might be effective in solving CKPs. We assume the strength of BR $c_n$
is fixed.

\medskip
\bth
\label{thm-cbrpt}
There exist functions $f_1, \dots, f_t$ with the following
property. Given $\svec \in \zad{t}$ with entries satisfying
 \beq
   \label{eq-sis}
   1 \leq s_t \leq \dots \leq s_1 \leq n-t,
 \eeq
if
 \beq
   \label{eq-Migtfi}
   M_i > f_i(M_{i+1}, \dots, M_t, s_i, P, \rvec, c_n),
   ~ i=1,\dots,t,
 \eeq
then
\beq
  \label{eq-ptzero}
  \tilde{P}_{i,1:s_i} = 0,~i=1,\dots,t.  \eeq Further, there exist
  $M_i$ with size polynomial in {\rm size}($P$), {\rm size}($\rvec$),
  and {\rm size}($c_n$) satisfying the Inequality {\rm
    (\ref{eq-Migtfi})}.  \enth

\medskip
Before presenting the somewhat technical proof of Theorem
\ref{thm-cbrpt}, we give some intuition for its implication. As an
example, consider the result in Equation (\ref{eq-ptzero}) for $n=12,
t=3$ (as in Example \ref{kp4}), and let $s_1=9, s_2 = s_3 = 8$. Then
the matrix $\tilde{P}$ has the form
\[
\begin{bmatrix}
  0 & 0 & 0 & 0 & 0 & 0 & 0 & 0 &  0  & \nz & \nz & \nz \\
  0 & 0 & 0 & 0 & 0 & 0 & 0 & 0 & \nz & \nz & \nz & \nz \\
  0 & 0 & 0 & 0 & 0 & 0 & 0 & 0 & \nz & \nz & \nz & \nz 
\end{bmatrix},~~\mbox{ where }~
\nz \mbox{ represent possible nonzero entries}.
\]
Intuitively, if $M_1$ is sufficiently larger than $M_2, \dots, M_t$,
then $\pvec_1 M_1$ contributes the most to the length of
$\avec$. Subsequently, if $M_2$ is sufficiently larger than $M_3,
\dots, M_t$, then the next biggest contribution to the norm of $\avec$
comes from $\pvec_2 M_2$, and so on. Hence, to shorten the columns of
$A$ in Equation (\ref{eq-Ackp}), the best option is to zero out
``many'' components of $\pvec_1$, followed by possibly fewer
components of $\pvec_2$, and so on. Since $P \xvec = P U \yvec =
\tilde{P} \yvec$, exploring all possible branches for $y_n, \dots,
y_{s_1+1}$ in the CBR reformulation is equivalent to branching on the
hyperplane defined by $\pvec_1 \xvec$ in the original CKP, exploring
all possible branches for $y_n, \dots, y_{s_2+1}$ is equivalent to
branching on the hyperplane defined by $\pvec_2 \xvec$, and so on.

We first present a lemma, which we use in the proof of Theorem
\ref{thm-cbrpt}. For brevity, we let
\begin{equation}
  \label{eq-alphak}
  \alpha_k = \Lambda_k(\NLt(P)),~k = 1,\dots,n-t.
\end{equation}
In words, $\alpha_k$ is the smallest number such that there are $k$
linearly independent vectors in $\NLt(P)$ with length bounded by
$\alpha_k$.

\ble
\label{lem-kthminLA}
For the matrix $A$ given in Equation {\rm (\ref{eq-Ackp})},
\beq
  \label{eq-kthminLA}
  \Lambda_k(\Lt(A)) \leq (\norm{\rvec} + 1) \alpha_k,~\mbox{ for }~
  k=1,\dots,n-t.
\eeq
\ele

\begin{proof}
 Let $k \leq n-t$, and let $\vvec_1, \dots, \vvec_k \in \NLt(P)$ be
 linearly independent vectors. Thus $\pvec_i \vvec_j = 0$ for all
 $i=1, \dots,t$ and $j=1,\dots,k$. Then $A \vvec_1, \dots, A \vvec_k
 \in \Lt(A)$ are linearly independent, and
 \[
 A \vvec_j = \begin{bmatrix} \avec \\ I \end{bmatrix} \vvec_j  =
 \begin{bmatrix}
  \pvec_1 M_1 + \pvec_2 M_2 + \dots + \pvec_t M_t + \rvec \\
  I \end{bmatrix} \vvec_j =
 \begin{bmatrix}
   \rvec \vvec_j \\ \vvec_j
 \end{bmatrix}~\forall j.
 \]
Thus we get $\norm{A\vvec_j} \leq (\norm{\rvec} + 1) \norm{\vvec_j} ~\forall
j$, and the bound in the Inequality (\ref{eq-kthminLA}) follows.
\end{proof}

\paragraph{Proof of Theorem \ref{thm-cbrpt}:} For brevity, we denote
$\rvec_i = \sum_{j=i}^t \pvec_j M_j + \rvec$ and $\tilde{\rvec}_i =
\rvec_i U$ for $i=1,\dots,t$. We also let $\rho = \max \{
\norm{\pvec_1}, \dots, \norm{\pvec_t},\, \norm{\rvec}+1\}$. We show
that
\[
f_i(M_{i+1}, \dots, M_t, s_i, P, \rvec, c_n) =
c_n \alpha_{s_i} \rho^2 (\,M_{i+1} + \dots + M_t + 1 \,)
\]
are suitable functions. Given that $\alpha_{s_i}$ have size polynomial
in the size of $P$ \cite{S86}, one could use these functions $f_i$ to
choose a set of $M_i$ that have size polynomial in the sizes of $P,
\rvec, c_n$.

Given $\svec \in \zad{t}$ satisfying the Inequalities (\ref{eq-sis}),
and assuming Inequality (\ref{eq-Migtfi}) holds, we show that
$\tilde{P}_{i,1:s_i} = 0$ holds for $i=1,\dots,\ell$ with $\ell \leq
t$ by induction, with the base case of $\ell = 0$ holding trivially by
extending the definitions to the case of $t=0$. Let $\ell \geq 1$, and
suppose this result holds for all $i < \ell$. We prove the result for
$i=\ell$.

Fix $j \leq s_{\ell}$. We are done if we manage to show
\beq
  \label{eq-ptlj0}
  \tilde{P}_{\ell, j} = 0.
\eeq
If $i < \ell$ and $j \leq s_{\ell}$, then $j \leq s_{\ell} \leq
s_i$. Hence the induction hypothesis implies $\tilde{P}_{i,j} =
0$. Recalling that $\tilde{a}_j$ denotes the $j$th entry of
$\tilde{\avec}$, we get
\[
\begin{aligned}
  \tilde{a}_j & = \sum_{i=1}^{\ell-1} \tilde{\pvec}_{i,j} M_i +
  \tilde{\pvec}_{\ell, j} M_{\ell} + \tilde{\rvec}_{\ell+1,j} \\
   & = \tilde{\pvec}_{\ell, j} M_{\ell} + \tilde{\rvec}_{\ell+1,j}.
\end{aligned}  
\]
To get a contradiction, assume Equation (\ref{eq-ptlj0}) does not
hold. Then we get
\[
\begin{aligned}
  \norm{\tilde{A}_{:,j} } \geq |\tilde{a}_j| & =
   |\tilde{\pvec}_{\ell, j} M_{\ell} + \tilde{\rvec}_{\ell+1,j}| \\
     & \geq
     |\tilde{\pvec}_{\ell, j} M_{\ell}| - |\tilde{\rvec}_{\ell+1,j}| \\
     & \geq M_{\ell} - |\tilde{\rvec}_{\ell+1,j}|,       
\end{aligned}  
\]
with the last inequality following from the fact that
$\tilde{\pvec}_{\ell, j} \in \Z$. Hence we get
\begin{equation}
  \label{eq-Mlub}
  \begin{aligned}
   M_{\ell}~ & \leq ~|\tilde{\rvec}_{\ell+1,j}| + \norm{\tilde{A}_{:,j}}\\
   & = ~ |\rvec_{\ell+1} U_{:,j}| + \norm{\tilde{A}_{:,j}} \\
   & \leq ~ \norm{\rvec_{\ell+1}} \,\norm{U_{:,j}} + \norm{\tilde{A}_{:,j}} \\
   & \leq ~ \norm{\rvec_{\ell+1}} \,\norm{\tilde{A}_{:,j}} + \norm{\tilde{A}_{:,j}} \\
   & = ~ (\norm{\rvec_{\ell+1}} + 1) \norm{\tilde{A}_{:,j}}.     
  \end{aligned}  
\end{equation}
The fourth inequality above, which replaced $\norm{U_{:,j}}$ with
$\norm{\tilde{A}_{:,j}}$, follows from the definition of $\tilde{A} =
AU$, which has $U$ as its submatrix (see Equation (\ref{eq-Ackp})).
Since $s_{\ell} \leq n-t$, there are $s_{\ell}$ linearly independent
vectors in $\NLt(P)$ with norm bounded by $\alpha_{s_{\ell}}$, and
hence by Lemma \ref{lem-kthminLA} there is the same number of linearly
independent vectors in $\Lt(A)$ with norm bounded by $(\norm{\rvec} +
1) \alpha_{s_{\ell}}$. Also, since $\tilde{A}$ was computed by BR with
strength $c_n$ and since $j \leq s_{\ell}$, we get that
\beq
  \label{eq-nrmAtj}
  \norm{\tilde{A}_{:,j}} \leq c_n (\norm{\rvec} + 1) \alpha_{s_{\ell}}.
\eeq
Combining the bounds in Inequalities (\ref{eq-Mlub}) and
(\ref{eq-nrmAtj}) yields
\[
  \begin{aligned}
    M_{\ell}~ & \leq ~ c_n \alpha_{s_{\ell}} (\norm{\rvec} + 1)
    (\,\norm{\rvec_{\ell+1}} + 1\,) \\
    & = ~ c_n \alpha_{s_{\ell}} (\norm{\rvec} + 1)
    (\,\norm{M_{\ell+1} \pvec_{\ell+1} + \dots + M_t \pvec_t + \rvec}
    + 1\,) \\
    & \leq ~ c_n \alpha_{s_{\ell}} (\norm{\rvec} + 1)
    \left(M_{\ell+1} \norm{\pvec_{\ell+1}} + \dots +
    M_t \norm{\pvec_t} + \norm{\rvec} + 1\,\right) \\
    & \leq ~ c_n \alpha_{s_{\ell}} \rho^2 \left(M_{\ell+1} + \dots +
    M_t  + 1\,\right), \\
  \end{aligned}  
\]
which provides the contradiction. \qed
\vspace*{-0.1in}

\section{Discussion} \label{sec_disc}

Restricting to directions defined by rational vectors, we could model
the problem of finding the direction along which the width of the
polyhedron of a given IP is the smallest as a mixed integer program---
see the work of Mahajan and Ralphs \cite{MR2010} for one such
model. At the same time, this MIP has more variables and constraints
than the original IP, and is typically harder to solve as well. For
instance, one could solve this MIP corresponding to the knapsack
instance (\ref{kp2}) to identify the knapsack coefficient vector
$\avec = (71,~82, \dots, 227)$ as the obvious thin direction, along
which the polyhedron has the minimal width of $0$. CPLEX takes $810$
\bb nodes to solve this MIP (as compared to $462$ nodes to solve the
original (CKP) itself---see Example \ref{ex2iw23ckp}). But we could
modify this MIP to identify another {\em useful} direction. For
$\pvec' = (26,~ 30,~ 34,~ 37,~ 49,~ 53,~ 56,~ 60,~ 72,~ 75,~ 79,~ 83)$
as identified by the modified MIP, we get $\width(\pvec',\ref{kp2}) =
358.98 - 358.66 = 0.32$ and hence $\iwidth(\pvec',\ref{kp2}) = 0$,
even though $\pvec' \xvec$ is not a thin direction.

CPLEX takes $2,644$ \bb nodes to solve the modified MIP that
identified $\pvec'$. Thus, if we {\em knew} $\pvec'$ beforehand, we
could solve (\ref{kp2}) at the root node by branching on the
hyperplane defined by $\pvec' \xvec$. At the same time, it appears
finding such a good branching direction is typically harder than
solving the original IP itself. Also notice that $\pvec'$ has larger
coefficients than $\pvec_1$ and $\pvec_2$. It would be interesting to
identify the class of IPs for which such a good direction is
guaranteed to exist. Also, would such a good direction have ``small''
coefficients?

We point out that our previous results on the hardness of ordinary \bb
on DKPs \cite{KrPa2009} (as well on more general integer-infeasible
knapsacks \cite{K07}) apply to the case of CKPs as well. Indirectly,
these results imply the hardness of \bb on CKPs when thin directions
as specified by the individual variables are used for branching.

While we demonstrated the effectiveness of CBR in solving the CKP
instances quickly, the main message we want to convey is the structure
of these problems: branching on the hyperplanes defined by
$\pvec_1\xvec, \pvec_2\xvec$, and $\pvec_3\xvec$ solves them quickly,
even though they might not be thin directions. On the other hand,
branching on thin directions (along the individual variables) might
not always be a good idea for \bbns. Indeed, if the structure of the
problem is assumed to be known, i.e., one is given $\pvec_1, \pvec_2,
\pvec_3, \rvec$ along with $M_1,M_2$, and $M_3$, one could verify
directly that branching on these hyperplanes solves the problem. But
if one is given just the final knapsack coefficient vector $\avec$,
CBR appears to be an effective method to {\em discover} that
structure. Alternatively, one could try to {\em guess} $\pvec_1$ from
$\avec$, e.g., using the ideas of diophantine approximation
\cite{S86}, add an extra variable that models $\pvec_1\xvec$, and
force CPLEX to branch on this extra variable. The idea of adding extra
variables in a similar setting was explored by Aardal and Wolsey
\cite{AaWo2010} in the context of lattice-based extended formulations
for integer inequality systems. Such an approach would also be much
more effective than trying to solve the original instances using
ordinary branch-and-cut.

\paragraph{{\large Acknowledgments:}} \normalsize The author thanks
G\'abor Pataki for useful discussions on the topics presented in this
paper. The author acknowledges partial support from the National
Science Foundation (NSF) through grant \#1064600.

\bibliographystyle{plain}

\end{document}